\documentclass{article}
\newcommand{\documentdate}{19 January 2025}

\usepackage{a4wide,latexsym,graphicx,amsmath}
\usepackage{varioref}
\usepackage{xcolor}
\usepackage{ulem}
\usepackage{hyperref}
\title{An optimally fast objective-function-free minimization
       algorithm using random subspaces}

\author{Stefania Bellavia\footnotemark[1],
        Serge Gratton\footnotemark[2],
        Benedetta Morini\footnotemark[1],
        Philippe L. Toint\footnotemark[3]}

\newcommand{\beqn}[1]{\begin{equation}\label{#1}}
\newcommand{\eeqn}{\end{equation}}
\newcommand{\req}[1]{(\ref{#1})}
\newcommand{\ms}{\;\;\;\;}
\newcommand{\tim}[1]{\;\; \mbox{#1} \;\;}

\newcounter{algo}[section]
\renewcommand{\thealgo}{\thesection.\arabic{algo}}
\newcommand{\llem}[2]{\vspace{\baselineskip} 
\noindent\framebox[\textwidth]{\parbox{0.95\textwidth}{
\begin{lemma} \label{#1} \rm #2 \end{lemma} } } \vspace{\baselineskip} }
\newcommand{\algo}[3]{\refstepcounter{algo}
\begin{center}\begin{figure}[htbp]
\framebox[\textwidth]{
\parbox{0.95\textwidth} {\vspace{\topsep}
{\bf Algorithm \thealgo : #2}\label{#1}\\
\vspace*{-\topsep} \mbox{ }\\
{#3} \vspace{\topsep} }}
\end{figure}\end{center}}
\newcommand{\bpr}{{\bf Proof.} \hspace{1.5mm}}
\newcommand{\epr}{\hfill $\Box$ \vspace*{1em}}
\newcommand{\lthm}[2]{\vspace{\baselineskip} 
\noindent\framebox[\textwidth]{\parbox{0.95\textwidth}{
\begin{theorem} \label{#1} \rm #2 \end{theorem} } } \vspace{\baselineskip} }

\newcommand{\mkB}{m_{k,B}}
\newcommand{\TB}{T_{k,B}}

\newcommand{\iiz}[1]{\{ 0, \ldots, #1 \}}
\newcommand{\iibe}[2]{\{ #1, \ldots, #2 \}}
\newcommand{\calA}{{\cal A}}

\newcommand{\calK}{{\cal K}} 
\newcommand{\calO}{{\cal O}} 
\newcommand{\calT}{{\cal T}} 
\newcommand{\calS}{{\cal S}} 
 
\newcommand{\calV}{{\cal V}}
\newcommand{\calW}{{\cal W}}

\newcommand{\calTa}{\calT_k^{(\omega)}}

\renewcommand{\Re}{\hbox{I\hskip -2pt R}}
\newcommand{\Na}{\hbox{I\hskip -2pt N}}
\newcommand{\prob}{\hbox{I\hskip -2pt P}}
\newcommand{\bigfrac}[2]{\frac{\displaystyle #1}{\displaystyle #2}}
\newcommand{\sfrac}[2]{{\scriptstyle \frac{#1}{#2}}}
\newcommand{\half}{\sfrac{1}{2}}
\newcommand{\eqdef}{\stackrel{\rm def}{=}}
\newcommand{\bigsum}{\displaystyle \sum}
\newcommand{\kap}[1]{\kappa_{\mbox{\tiny #1}}}
\newcommand{\khigh}{\kappa_{\mbox{\tiny high}}}
\newcommand{\al}[1]{{\footnotesize{\sf #1}}}
\newcommand{\tal}[1]{{\normalsize {\sf #1}}}

\newcommand{\pione}{\pi_S^{(1)}}
\newcommand{\pioneone}{\pi_1^{(1)}}

\newcommand{\ree}[1]{\overset{\req{#1}}{=}}
\newcommand{\rge}[1]{\overset{\req{#1}}{\ge}}
\newcommand{\rgt}[1]{\overset{\req{#1}}{>}}
\newcommand{\rle}[1]{\overset{\req{#1}}{\le}}

\newcommand{\skoffar}{\al{SKOFFAR}}
\newcommand{\skoffarp}{\al{SKOFFAR$p$}}

\newcommand{\skoffarduedot}{\al{\skoffar{\tiny{$2$}}}}
\newcommand{\skoffardue}{\al{\skoffar{\tiny{$2$}}} }
\newcommand{\tSKOFFARp}{\tal{SKOFFAR$p$}}

\newtheorem{theorem}{Theorem}[section]
\newtheorem{lemma}[theorem]{Lemma}

\newtheorem{corollary}[theorem]{Corollary}
\newtheorem{definition}[theorem]{Definition}

\newcommand{\proof}[1]{
\begin{list}{}{
\setlength{\topsep}{0.0pt}
\setlength{\partopsep}{0.0pt}
\setlength{\leftmargin}{0.025\textwidth}
\setlength{\rightmargin}{0.5\leftmargin}
\setlength{\labelwidth}{0.5\leftmargin}
\setlength{\labelsep}{0.25\leftmargin}}
\item \bpr #1 \epr \noindent
\end{list}}

\begin{document}

\maketitle
\centerline{\documentdate}

\footnotetext[1]{Dipartimento di Ingegneria Industriale,
    Universit\`{a} degli Studi di Firenze, Italy. Member of the INdAM Research
    Group GNCS. Email: stefania.bellavia@unifi.it, benedetta.morini@unifi.it.
    Work partially supported by INDAM-GNCS through Progetti di Ricerca 2024 and 
by PNRR - Missione 4 Istruzione e Ricerca - Componente C2 Investimento 1.1, Fondo per il Programma Nazionale di Ricerca e Progetti di Rilevante Interesse Nazionale (PRIN) funded by the European Commission under the NextGeneration EU programme, project 
``Numerical Optimization with Adaptive Accuracy and Applications to Machine Learning'',  code: 2022N3ZNAX
 MUR D.D. financing decree n. 973 of 30th June 2023 (CUP B53D23012670006), and by Partenariato esteso FAIR ``Future Artificial Intelligence Research'' SPOKE 1 Human-Centered AI. Obiettivo 4, Project ``Mathematical and Physical approaches to innovative Machine Learning technologies (MaPLe)''.}
\footnotetext[2]{Universit\'e de Toulouse, INP, IRIT, Toulouse,
  France. Work partially supported by 3IA Artificial and Natural
  Intelligence Toulouse Institute (ANITI), French ``Investing for the
  Future - PIA3'' program under the Grant agreement ANR-19-PI3A-0004.
  Email: serge.gratton@enseeiht.fr.}
\footnotetext[3]{Namur Center for Complex Systems (naXys),
  University of Namur, 61, rue de Bruxelles, B-5000 Namur, Belgium.
  Email: philippe.toint@unamur.be.}

\begin{abstract}
An algorithm for unconstrained non-convex optimization is described, which does not evaluate the objective function and in which minimization is carried out, at each iteration, within a randomly selected subspace. It is shown that this random approximation technique does not affect the method's convergence nor its evaluation complexity for the search of an $\epsilon$-approximate first-order critical point, which is $\mathcal{O}(\epsilon^{-(p+1)/p})$, where $p$ is the order of derivatives used. A variant of the algorithm using approximate Hessian matrices
is also analysed and shown to require at most $\mathcal{O}(\epsilon^{-2})$ evaluations.
Preliminary numerical tests show that the random-subspace technique can significantly improve performance when used with $p=2$ in the correct context, making it very competitive when compared to standard first-order algorithms.
\end{abstract}

{\small
\textbf{Keywords:} nonlinear optimization, stochastic adaptive 
regularisation methods, sketching, evaluation complexity, 
objective-function-free optimization (OFFO).
}

\section{Introduction}

Recent years have seen the emergence of random concepts in iterative
algorithms for nonconvex optimization (see \cite{CS2020} and reference therein and 
\cite{CBS2022,BGMT2022,BKMR2023,BW2023,SX2023}).
In particular, several authors 
\cite{VuPoirDambLibe19,FRW2021,Shao21,DW2022,CartFowkShao22,CartShaoTans25}\footnote{\cite{CartShaoTans25} was posted on arXiv when the revision of the present
paper was being finalized. It elaborates on the results of 
\cite{Shao21,CartFowkShao22} by using a variable dimension sketching strategy.}  have
suggested algorithms in which the search for a better iterate is carried
out in random subspaces of the space of variables, instead of, as is
more traditional and often more costly, in the complete space. In these
proposals, the Johnson-Lindenstrauss embedding Lemma
(see \cite{DasgGupt02} for a simple exposition)
is used to ensure that the relevant information can be found very
efficiently in the selected subspace with high probability, and this leads to an elegant analysis yielding optimal complexity bounds
for ``random-subspace'' variants of the standard trust-region and adaptive-regularisation methods for unconstrained minimization. In parallel
with this interesting development, alternative non-standard
optimization methods have also been introduced where the
objective function of the problem is never computed (these algorithms
use derivatives' values only). The motivation for such methods 
originates in applications with noisy objective functions. Indeed, because 
differences of objective function's values are not used to accept 
or reject iterates, the methods' behaviour is much less sensitive to noise 
than that of the more standard algorithms using function 
values \cite{GratJeraToin23c}.
This new class of ``objective-function-free optimization'' (OFFO) methods 
includes popular first-order algorithms
as \al{ADAM} or \al{ADAGRAD}, and has been investigated, for instance, in 
\cite{DuchHazaSing11,KingBa15,TielHint12,GrapStel22,WuWardBott18}.

The purpose of this paper is to discuss an
algorithm which combines these two ideas for the first time
while maintaining the desirable properties of both.  More specifically, we
describe an OFFO adaptive regularisation method using first- or
higher-order models defined in random subspaces, and show that this
algorithm still enjoys the optimal global rate of convergence known
for comparable adaptive-regularisation methods.  Independently of the
practical interest for such a method,
which, we argue below, can be substantial in the right context,
our analysis is a new step in
the ``information thinning'' question, which is to isolate what
information is necessary for a minimization method to achieve optimal
complexity.  Indeed, while \cite{GratJeraToin22c} proves that function
values are unnecessary, the present paper further shows that this is
also the case for ``full space'' information\footnote{One might argue
that it has long been known that information along the directions
given by the gradient and the step suffices, but this requires the
step to be known and thus amounts to an \textit{a posteriori} observation
instead of an \textit{a priori} algorithmically exploitable strategy.}
under suitable probabilistic assumptions.

Our approach has a further advantage compared to existing proposals,
like the random-subspace trust-region and random-subspace regularisation methods
of \cite{Shao21} and \cite{CartFowkShao22}.
Because no evaluation of the objective function is involved, the algorithm
generates a much simpler random process (there is now only one random
event per iteration), in turn considerably simplifying the proofs as
the number of iteration types whose number must be estimated
(in \cite[Chapter 4]{Shao21}) is now reduced to only two.
While our theory covers the general case where derivatives of higher order 
than one are estimated, our practical focus will be on  the case where  
first and second derivatives are used.

The paper is organised as follows. The new algorithm is proposed in
Section~\ref{algo-s}, while its evaluation complexity is analysed
under general embedding conditions in Section~\ref{complexity-s}.
A brief discussion of a possible way to select the random subspaces
are presented in Section~\ref{select-s}.
The numerical behaviour of the second-order variant is illustrated 
in Section~\ref{numerics-s}. Some conclusions are finally presented 
in Section~\ref{concl-s}. A discussion of a further variant using quadratically 
regularised inexact quadratic models is proposed and analysed in appendix.

\section{An OFFO adaptive regularisation algorithm using random subspaces}\label{algo-s}

\noindent
The problem of interest in what follows is the standard nonconvex
unconstrained minimization of a (sufficiently) smooth objective
function, that is 
\beqn{problem}
\min_{x\in {\scriptsize{\Re^n}}}f(x),
\eeqn
where $f:\Re^n \rightarrow \Re$.
As indicated in the introduction, our aim is to design an adaptive
regularisation algorithm in which \textit{the objective function value is
never computed}, and in which the step is obtained by approximately
minimizing a suitable model of the objective function in a random
subspace. To ensure that this approach is sensible, we make the
following assumptions. 

\vspace*{5mm}
\noindent
\textbf{AS.1} $f$ is $p$ times continuously differentiable in $\Re^n$. 

\noindent
\textbf{AS.2} There exists a constant $f_{\rm low}$ such that
$f(x) \geq f_{\rm low}$ for all $x\in \Re^n$.

\noindent
\textbf{AS.3} The $p$th derivative of $f$ is globally Lipschitz
continuous, that is, there exists a non-negative constant $L_p$ such that
\[
\|\nabla_x^p f(x) - \nabla_x^p f(y)\| \leq L_p \|x-y\| \, \text{ for all } x,y \in \Re^n,
\]
where $\|.\| $ denotes the Euclidean norm for vectors
in $\Re^n$ and the corresponding subordinate norm for tensors. 

\noindent
\textbf{AS.4}
The gradient of $f$ is bounded, that is there exists a constant
$\kap{g}\geq 0$ such that, for all $x\in\Re^n$,
\[
\|\nabla_x^1f(x)\| \leq \kap{g}.
\]

\noindent
\textbf{AS.5}
If $p > 1$, there exists a constant $\khigh\geq 0$ such that
\[
\min_{\|d\|\leq1} \nabla_x^i f(x) [d]^i \geq -\khigh \tim{for all}
x\in\Re^n \tim{and} i\in\iibe{2}{p},
\]
where $\nabla_x^i f(x)$ is the $i$th derivative tensor of
$f$ computed at $x$, and where $T[d]^i$ denotes the
$i$-dimensional tensor $T$ applied on $i$ copies of the vector $d$.
(For notational convenience, we set $\khigh = 0$ if $p=1$.)

\vspace*{5mm}
\noindent
We refer the reader to
\cite[Appendix~6]{CartGoulToin22} for details on derivative tensors.
Observe that, given AS.1, AS.3 is automatically satisfied if the iterates generated by the algorithms remain in a bounded domain. This is in particular the case if a objective function's level set identified in Lemma~\ref{fk1boundsigmak1} is bounded because, as we comment on after this lemma, it contains all generated iterates.
Observe also that AS.5 is irrelevant in the case where $p=1$.
Should one be interested in higher-order methods, 
AS.5 is weaker than assuming uniform boundedness of the
derivative tensors of degree two and above (there is no upper bound on
the value of $\nabla_x^if(x) [d]^i$), or, equivalently, Lipschitz
continuity of derivatives of degree one to $p-1$.

\subsection{The \tSKOFFARp\ algorithm}

As suggested above,
adaptive regularisation methods are iterative schemes which compute a
step from an iterate $x_k$ to the next by approximately minimizing a
$p$-th degree regularised model $m_k(s)$ of $f(x_k+s)$ of the form
\beqn{model}
m_k(s) \eqdef T_{f,p}(x_k,s) + \frac{\sigma_k}{(p+1)!} \|s\|^{p+1},
\eeqn
where $T_{f,p}(x,s)$ is the $p$th order Taylor expansion of
functional $f$ at $x$ truncated at order $p$, that is,
\beqn{taylor model}
T_{f,p}(x,s) \eqdef f(x) + \sum_{i=1}^{p} \frac{1}{i!}\nabla_x^i f(x)[s]^i.
\eeqn
To obtain the model \req{model}, the $p$-th order Taylor series
\req{taylor model} is ``regularised'' by adding the term
$\frac{\sigma_k}{(p+1)!} \|s\|^{p+1}$ (where $\sigma_k$ is the
iteration-dependent regularisation parameter), thereby ensuring that
$m_k(s)$ is bounded below and that a step $s_k$ (approximately)
minimizing this model is well-defined.  

Following \cite{Shao21}, we propose to compute a random subspace step
at iteration $k$ as follows. Given an iteration-independent
distribution $\calS$ of $\ell \times n$ random matrices (with
$\ell<n$), let $S_k$ be drawn from this distribution and consider
minimizing the sketched regularised model 
\beqn{model sub}
\widehat m_k(\widehat s) \eqdef \widehat T_{f,p}(x_k,\widehat s) + \frac{\sigma_k}{(p+1)!} \|S_k^T \widehat s\|^{p+1},
\eeqn
as a function of $\widehat s\in \Re^{\ell}$, where the 
sketched Taylor model $\widehat T_{f,p}(x,\widehat s)$ is given by
\[
\widehat T_{f,p}(x,\widehat s) \eqdef f(x) + \sum_{i=1}^{p} \frac{1}{i!}\nabla_x^i f(x)[S_k^T \widehat s]^i.
\]
Letting $\widehat s_k$ an approximate minimizer of $\widehat m_k(\widehat s)$, the full dimensional step is then defined by  $s_k=S_k^T \widehat s_k$.
We note that $\widehat T_{f,p}(x_k,\widehat s_k)= T_{f,p}(x_k, s_k)$ and 
\beqn{equalmodels}
\widehat m_k(\widehat s_k)= m_k( s_k).
\eeqn

\algo{SKOFFARp}{Sketching OFFO adaptive regularisation of degree $p$ (\skoffarp)}{
\begin{description}
\item[Step 0: Initialization: ] An initial point $x_0\in \Re^n$, a regularisation
parameter $\nu_0>0$ and a requested final gradient accuracy
$\epsilon \in (0,1]$ are given, as well as the parameters
\[
\theta > 1, \;\mu_{-1} \geq 0 \tim{ and }   0< \vartheta <1. 
\]
Set $k=0$.
\item[Step 1: Step calculation: ]
If $k=0$, set $\sigma_0=\nu_0$.
Otherwise, select
\beqn{sigkupdate}
\sigma_k  \in \Big[ \vartheta \nu_k, \max[\nu_k, \mu_k] \Big],
\eeqn
where
\beqn{muk-def}
\mu_k = \max\left[ \mu_{k-1}, \frac{\|S_{k-1}g_k\| - \|\nabla^1_{\widehat s}\widehat T_{f,p}(x_{k-1},\widehat s_{k-1})\|}{\kappa_{S,k-1}.\|s_{k-1}\|^p}\right],
\eeqn
with some $\kappa_{S,k-1}$ such that $\|S_{k-1}\| \leq \kappa_{S,k-1}$.
Draw a random matrix $S_k\in \Re^{\ell \times n}$ from
$\calS$ and compute a step $s_k=S_k^T \widehat s_k$  such
that $\widehat s_k$ sufficiently reduces  the model 
$\widehat m_k$ defined in \req{model sub} in
the sense that  
\beqn{descent}
 \widehat m_k(\widehat s_k) - \widehat m_k(0) <  0
\eeqn
and
\beqn{gradstep}
\|\nabla_{\widehat s}^1 \widehat T_{f,p}(x_k,\widehat s_k)\|
\leq  \theta \frac{\sigma_k}{p!}\|S_k^T\widehat s_k\|^{p-1}\|S_kS_k^T\widehat s_k\|. 
\eeqn
\item[Step 2: Updates. ]
Set
\[
x_{k+1} = x_k + s_k
\]
and
\beqn{vkupdate}
\nu_{k+1} =  \nu_k + \nu_k \| s_k\|^{p+1}.
\eeqn
Increment $k$ by one and go to Step~1.
\end{description}
}

Some comments on this algorithm are necessary.
\begin{enumerate}
\item It is crucial to observe that, while the definition of the model
  in \req{model sub} involves the function value $f(x_k)$ (in
  $\widehat T_{f,p}(x_k,\widehat s)$), this function value is never
  needed in the algorithm (it cancels out in \req{descent}) and
  therefore need not to be evaluated. The algorithm thus belong to the
  OFFO class. Of course, the minimization of the model may require the
  evaluation of the sketched derivatives $\{\nabla_{\widehat s}^j
  f(x_k)[S_k\cdot]^j\}_{j=1}^p$, at least along some
  directions\footnote{In the course of a Krylov subproblem solver for
  $p=2$, say.}. This makes the use of derivatives of degree higher
  than two potentially useable in practice, especially if the
  objective function is partially separable \cite{GrieToin82b,ConnGoulToin92}. 
\item
  Since 
  \[
  \nabla^1_{\widehat s} \|S_k^T \widehat s\|^{p+1}
  = (p+1) \|S_k^T \widehat s\|^{p-1}\,S_kS_k^T \widehat s,
  \]
  one verifies that conditions \req{descent} and \req{gradstep} do hold at an
  exact minimizer of $\widehat m_k$ (the latter with $\theta=1$). A
  step satisfying these conditions is therefore guaranteed to exist.
  Note that \req{gradstep} is a condition on the norm of the gradient 
  of the Taylor series for $f$, at variance with \cite{CartGoulToin22,Shao21}
  where the condition is on the gradient of the regularised model \req{model}.
\item
  At variance with standard trust-region and adaptive-regularisation
  methods, the algorithm does not involve any (typically noise
  sensitive) test to accept or reject the trial iterate $x_k+s_k$, and
  every trial point is thus ``accepted'' as the new iterate. In the
  vocabulary used for trust-region and adaptive regularisation
  methods, every iteration is therefore ``successful''. 
\item 
  The value of $\mu_k$ in the definition \req{sigkupdate} of $\sigma_k$ 
  is chosen to help the regularisation parameter 
  $\sigma_k$ to grow fast enough, given the knowledge at iteration $k$.
  We will show in 
  Lemma~\ref{muk-bound} that $\mu_k$ is bounded above by $\max[\mu_{-1},L_p]$ irrespective of the choice of $\kappa_{S,k-1}$.
As a consequence, the specific values of $\kappa_{S,k-1}$  in \eqref{muk-def}  play no role in our 
complexity analysis, albeit they obviously affect the practical performance
of the method. Finally, we stress that the knowledge of the constants
$L_p$ and $\kappa_g$, given in AS.3 and AS.4 respectively, is
\textit{not required} in the algorithm.
\end{enumerate}

The \skoffarp\ algorithm can be seen as a stochastic
process because the selection of $S_k$ is
 random and  yields random
realizations\footnote{Formally, the iterates and steps are
random variables on some implicitly defined probability space, and
$x_k$ and $s_k$ are their realizations.} of the iterates $x_k$ and of the
steps $s_k$. The objective of our forthcoming complexity analysis for
this algorithm is to derive a probabilistic bound on the process
hitting time
\beqn{kepsdef}
N_1(\epsilon) \eqdef \min\{ k \in \Na \mid \|g_k\| \leq \epsilon \},
\eeqn
where we denote $g_k \eqdef \nabla_x^1 f(x_k)$ for all $k$.
 $N_1(\epsilon)$ is the number of iterations that a particular
realization of the algorithm requires to obtain an
$\epsilon$-approximate first-order critical point. 

\section{Evaluation complexity for the \tSKOFFARp\ algorithm}\label{complexity-s}

\noindent
Before discussing our analysis of evaluation complexity, we first
restate some classical lemmas for \al{AR$p$} algorithms, starting with
Lipschitz error bounds. 

\llem{lipschitz}{
  Suppose that AS.1 and AS.3 hold. Then
  \beqn{Lip-f}
  f(x_{k+1})- \widehat T_{f,p}(x_k,\widehat s_k)
  = f(x_{k+1})- T_{f,p}(x_k,s_k) \leq \frac{L_p}{(p+1)!} \|s_k\|^{p+1},
  \eeqn
  and
  \beqn{Lip-g}
  \|g_{k+1}-\nabla_s^1 T_{f,p}(x_k,s_k) \| \leq \frac{L_p}{p!} \|s_k\|^p.
  \eeqn
}

\proof{This is a standard result (see \cite[Lemma~2.1]{CartGoulToin20b}
  for instance).}

\noindent
We next state a simple lower bound on the Taylor series' decrease.

\llem{model-decrease}{
  \beqn{Tdecr}
  \Delta T_{f,p}(x_k,s_k) \eqdef T_{f,p}(x_k,0)-T_{f,p}(x_k,s_k)
  > \frac{\sigma_k}{(p+1)!} \|s_k\|^{p+1}.
  \eeqn
}

\proof{The bound directly results from 
$\widehat m_k(\widehat s_k)= m_k(s_k)$, \req{descent} and \req{model}.}

\noindent
This and AS.2 allow us to establish a lower bound on the decrease in
the objective function (although it is never computed).

\llem{Ifsigmabig}{
Suppose that AS.1 and AS.3 hold and that $\sigma_k \geq 2 L_p$. Then
\beqn{sigma-upper}
f(x_{k}) - f(x_{k+1}) > \frac{\sigma_k}{2(p+1)!} \|s_k\|^{p+1}.
\eeqn
}

\proof{
  From \req{Lip-f} and \req{Tdecr}, we obtain that
  \[
  f(x_{k}) - f(x_{k+1}) >
  \frac{\sigma_k-L_p}{(p+1)!}\|s_k\|^{p+1}
  \]
  and \req{sigma-upper} immediately follows from our assumption on $\sigma_k$.
} 

\noindent
We now recall an upper bound on $\|s_k\|$ generalizing those proposed
in \cite{CartGoulToin11,GratToin21} to the case where $p$ is arbitrary.

\llem{stepgkbound}{
Suppose that AS.1 and AS.5 hold. At each iteration $k$, we have that
\beqn{skbound}
\| s_k\| \leq
2 \eta +2 \left(\frac{(p+1)! \|g_k \|}{\sigma_k}\right)^\sfrac{1}{p},
\eeqn
where 
\beqn{kaphigh}
\eta = \sum_{i=2}^{p}  \left[\frac{\khigh(p+1)!}{i!\,\vartheta \nu_0}\right]^\sfrac{1}{p-i+1}.
\eeqn
}
\proof{See \cite[Lemma 3.6]{GratJeraToin22c}. Note that this result does
not involve $S_k$ as it is valid for any step which 
reduces $m_k$ and, using 
\req{equalmodels} and \req{descent}, 
$m_k(s_k)=\widehat{m}_k(\widehat{s}_k)< \widehat{m}_k(0) = m_k(0)$.
}

\noindent
Our next step is to show that $\mu_k$ is bounded.

\llem{muk-bound}{Suppose that AS.1 and AS.3 hold.
For all $k\geq 0$, 
\beqn{muk-bounded}
\mu_k \leq \max[\mu_{-1},L_p ].
\eeqn
}
\proof{
We have that 
$\nabla_{\widehat{s}}^1\widehat{T}_{f,p}(x_{k-1},\widehat{s}_{k-1})
=\nabla_{\widehat{s}}^1T_{f,p}(x_{k-1},S_{k-1}^Ts_{k-1})
=S_{k-1}\nabla_s^1 T_{f,p}(x_{k-1},s_{k-1})$, 
so that, using the triangular inequality, \req{Lip-g} and \req{gradstep},
\[
\begin{array}{lcl}
\|S_{k-1}g_k\|
& \leq & \|S_{k-1} (g_k -\nabla_x^1 T_{f,p}(x_{k-1},s_{k-1}))\|
   + \|S_{k-1}\nabla_x^1 T_{f,p}(x_{k-1},s_{k-1})\|\\*[1ex]
& \leq & \|S_{k-1}\| L_p \|s_{k-1}\|^p 
    + \| \nabla_{\widehat{s}}^1\widehat{T}_{f,p}(x_{k-1},\widehat{s}_{k-1})\|,
\end{array}
\]
and thus
\beqn{Lpmuk}
L_p 
\geq \frac{\|S_{k-1}g_k\|-\| \nabla_{\widehat{s}}^1\widehat{T}_{f,p}(x_{k-1},\widehat{s}_{k-1})\|}{\|S_{k-1}\|\|s_{k-1}\|^p}.
\eeqn
The inequality \req{muk-bounded} then follows from \req{muk-def} and 
$\|S_{k-1}\|\leq\kappa_{S,k-1}$.
}

\noindent
The proof of this lemma  shows that a tighter lower bound 
on $L_p$ (see \req{Lpmuk}) is also available at the often 
significant cost of evaluating $\|S_{k-1}\|$, thus motivating the 
introduction of the (hopefully) cheaper $\kappa_{S,k-1}$.

Since our objective is to minimize $f$, obtaining a decrease as stated
by Lemma~\ref{Ifsigmabig} is important.  The condition $\sigma_k \geq 2 L_p$ in this lemma and \req{sigkupdate} together suggest that the condition
\beqn{vkbig}
\nu_k \geq \frac{2 L_p}{\vartheta}
\eeqn
is important for our subsequent analysis. Remembering that
$\nu_k$ is increasing with $k$, we therefore define
\beqn{k1def}
k_1 \eqdef \inf\left\{ k\geq 1 \mid \nu_k \geq \frac{2L_p}{\vartheta}\right\}
\eeqn
the index of the first iterate (in a given realization) such that significant objective function decrease is
guaranteed by Lemma~\ref{Ifsigmabig}.  Note that $k_1$ may be infinite, which is why we define the random event
\beqn{calK1def}
\calK_1 \eqdef \{ k_1 \; \mbox{as defined by \req{k1def} is  finite}\}.
\eeqn
We now pursue our analysis under the condition that $\calK_1$ occurs.
The next series of lemmas provides bounds, conditional on $\calK_1$,  on $f(x_{k_1})$ and
$\nu_{k_1}$,  which in turn allows establishing an upper bound on
the regularisation parameter, only depending on the problem and the
fixed algorithmic parameters.

\llem{vk1bound}{
  Suppose that AS.1, AS.3, AS.4 and AS.5 hold
  and consider a realization of the \skoffarp\ algorithm where
  $\calK_1$ occurs.
  Then
\beqn{numaxdef}
 \nu_{k_1} \leq \nu_{\max}
\eqdef 
 \frac{2L_p}{\vartheta}\left[1+\left(2 \eta+2\left(\frac{(p+1)!\kap{g}}{\vartheta \nu_0}\right)^\sfrac{1}{p}\right)^{p+1}\right], 
\eeqn
where $\eta$ is defined in \req{kaphigh} and $\kap{g}$ in AS.4.
}

\proof{
Since $\calK_1$ is assumed to occur, $k_1$ is well-defined and finite.
Successively using Lemma~\ref{stepgkbound} and the update rule for
$\nu_k$ \req{vkupdate}, we derive that
\[
\nu_{k_1} \ree{vkupdate} \nu_{k_1 - 1} + \nu_{k_1 -1} \|s_{k_1 -1} \|^{p+1} 
\rle{skbound}  \nu_{k_1 - 1} +  \nu_{k_1 -1}
  \left( 2 \left( (p+1)! \frac{\| g_{k_1-1}\|}{\sigma_{k_1-1}}\right)^\sfrac{1}{p} + 2\eta\right)^{p+1}
  \]
and the desired result follows by using AS.4, the definition of
$k_1$ in \req{k1def} and the inequalities $\sigma_{k_1-1} \ge \vartheta \nu_{k_1-1}\ge \vartheta \nu_0$.
}

\noindent
Lemma~\ref{vk1bound} allows us to establish an upper bound on $f(x_{k_1})$ as a function of $\nu_{\max}$.

\llem{fk1boundsigmak1}{
  Suppose that AS.1, AS.3, AS.4 and AS.5 hold and consider a
  realization of the \skoffarp\ algorithm where $\calK_1$ occurs. 
  Then
\beqn{fmaxdef}
f(x_{k_1}) \leq f_{\max} \eqdef f(x_0) +
          \frac{1}{(p+1)!}\left(\frac{L_p}{\sigma_0} \nu_{\max}
          +\vartheta \sigma_0\right).
\eeqn
}
\proof{
Lemma 3.8 in \cite {GratJeraToin22c}  shows that, for any $k\geq 0$,
\[
f(x_k) \leq  f(x_0) + \frac{1}{(p+1)!}\left( \frac{L_p\nu_k}{\sigma_0} +\vartheta\sigma_0)\right).
\]
The desired bound then follows from Lemma \ref{vk1bound}.
}

\noindent
Observe that this result ensures that all iterates generated by the algorithm belong to the level set $\{x\in \Re^n \mid f(x) \leq f(x_{k_1})\}$.
The two bounds stated in Lemmas~\ref{fk1boundsigmak1} and
\ref{vk1bound} are also useful in that they now imply an upper bound on
the regularisation parameter, an important step in standard theory for
regularisation methods.

\llem{uppersigmak}{
Suppose that AS.1, AS.2, AS.3, AS.4 and AS.5 hold and consider a 
realization of the \skoffarp\ algorithm. Then
\beqn{sigmamaxdef}
\begin{array}{lcl}
\sigma_k & \leq  &\sigma_{\max}\\
& \eqdef  &\max\left[
\bigfrac{2(p+1)!}{\vartheta}\left[f(x_0)-f_{\rm low}
+\bigfrac{1}{(p+1)!}\left(\bigfrac{L_p}{\sigma_0}
\nu_{\max}+\vartheta\sigma_0 \right)\right]+ \nu_{\max},
\mu_{-1}, \bigfrac{2 L_p}{\vartheta}, \nu_0 \right].
\end{array}
\eeqn
}

\proof{We proceed as in Lemma \cite[Lemma 3.9]{GratJeraToin22c} and give the proof for sake of clarity. Suppose first that $\calK_1$ occurs. From the definition of $k_1$ in \eqref{k1def}, we deduce that $\sigma_j \geq 2L_p$.  From Lemma~\ref{Ifsigmabig}, we then have that  
\[
f(x_{j}) - f(x_{j+1})
\geq \frac{\sigma_j}{2(p+1)!} \|s_j \|^{p+1}
\geq \vartheta  \frac{\nu_j}{2(p+1)!} \|s_j \|^{p+1}.
\] 
Summing the previous inequality from $j = k_1$ to $k-1$ and using the
$\nu_j$ update rule  \eqref{vkupdate} and AS.2, we deduce that
\[
f(x_{k_1}) - f_{\rm low}\geq f(x_{k_1}) - f(x_k) \geq  \frac{\vartheta}{2(p+1)!} (\nu_k - \nu_{k_1}).
\]
Rearranging the previous inequality and using Lemma~\ref{vk1bound} then gives that
\beqn{nuk-upper}
\nu_k \leq \frac{2(p+1)!}{\vartheta} \left( f(x_{k_1}) - f_{\rm low} \right) + \nu_{\max}.
\eeqn
Combining now Lemma~\ref{fk1boundsigmak1} (to bound $f(x_{k_1})$),
\req{sigkupdate} and \req{muk-bounded} yields that
\[
\sigma_k \leq  \sigma_{\max}
\eqdef \max\left[
\bigfrac{2(p+1)!}{\vartheta}\left[f(x_0)-f_{\rm
  low}+\frac{1}{(p+1)!}\left(\frac{L_p}{\sigma_0}
\nu_{\max}+\vartheta\sigma_0 \right)\right]+ \nu_{\max},
\mu_{-1},  L_p,  \nu_0 \right].
\]
If $\calK_1$ does not occur,
$\nu_k\le 2L_p/\vartheta$ for all $k$. Thus we obtain,
using  \req{sigkupdate} and \req{muk-bounded}, that 
$\sigma_k\le \max[\frac{2L_p}{\theta},\mu_{-1}]$ for all $k$,
and \req{sigmamaxdef} also holds.
}

\noindent
The theory of adaptive-regularisation methods crucially depends on the
relation between the steplength $\|s_k\|$ and the norm of the gradient
at the next iteration $\|g_{k+1}\|$ (see Lemmas~3.3.3 and 4.1.3 in
\cite{CartGoulToin22}, for instance), which is itself bounded 
below by $\epsilon$ before convergence.  Here we choose to consider this
dependence as a random event, depending on the choice of $S_k$.
This is formalized in the following definition.

\begin{definition}\label{deftrue}
Given some $\epsilon \in (0,1)$ independent of $k$,
iteration $k\in \iiz{N_1(\epsilon)-2}$ is said to be $\omega$-true for some $\omega
\in (0,1)$ independent of $k$ whenever
\beqn{crucial}
\|s_k\|^p \geq \omega \epsilon.
\eeqn
\end{definition}
\vskip 5pt

\noindent
We discuss in Section~\ref{select-s} conditions which may enforce this
property, but immediately note that it automatically holds
if $S_k$ is of rank $n$ \cite[Lemma~3.4]{GratJeraToin22c}.
We also define
\beqn{set-true}
\calTa\eqdef \{ j\in \iiz{k-1} \mid \mbox{iteration $j$ is  $\omega$-true}\},
\eeqn
the index set $\calTa$ of all $\omega$-true iterations in the first $k$.

\noindent

Given these definitions, we now need to establish under which condition the event $\calK_1$ occurs with
high probability. Such a condition is obtained in two stages, the
first follows arguments by \cite[Lemma~7]{GrapStel22} and
\cite[Lemma~3.5]{GratJeraToin22c} and investigates, in our probabilistic
setting, the effect of accumulating $\omega$-true iterations.

\llem{sigmanotbig}{
Suppose that AS.1 and AS.3 hold and consider a particular realization
of the \skoffarp\ algorithm.  Let $k_0< N_1(\epsilon)$ be an
iteration index (in this realization) such that  $k_*$ $\omega$-true
iterations have been performed among those of index 0 to $k_0-1$,
where 
\beqn{kstardef}
k_*
\eqdef
\left\lceil\frac{2L_p\epsilon^{-\sfrac{p+1}{p}}}{\vartheta \nu_0\, \omega^{\sfrac{p+1}{p}}}
\right\rceil.
\eeqn
Then $k_1$ exists, $k_1\leq k_0$
and, for all $k\geq k_1$, 
\beqn{sigmabig}
\sigma_k \geq 2L_p.
\eeqn
\vspace*{-5mm}
}
\proof{
First observe that \req{sigmabig} is a direct consequence of
\req{sigkupdate} if $\nu_k \geq 2L_p/\vartheta$. Suppose now that,
for some $k \in \iibe{k_0}{N_1(\epsilon)-1}$, $\nu_k < 2L_p/\vartheta$. 
Since $\{\nu_k\}$ is a
non-decreasing sequence, we deduce that this inequality holds for
$j \in\iiz{k}$. Successively using the form of the $\nu_k$ update
rule \req{vkupdate}, \req{crucial}, \req{sigkupdate} and the fact
that $k<N_1(\epsilon)$, we obtain that
\begin{align*}
\nu_k
 & \rgt{vkupdate} \sum_{j=0}^{k-1} \nu_j \|s_j\|^{p+1}
   \rgt{set-true} \sum_{j\in \calTa} \nu_j \| s_j\|^{p+1}
   \rge{crucial}  \sum_{j\in \calTa} \nu_j  (\omega\epsilon)^{\sfrac{p+1}{p}}\\
 & \rge{sigkupdate} \sum_{j\in \calTa} \nu_0 \left(\omega \epsilon\right)^{\sfrac{p+1}{p}}
  \\rge{set-true} k_*\, \nu_0 (\omega \epsilon)^{\sfrac{p+1}{p}}.
\end{align*}
Substituting the definition of $k_*$ in the last inequality, we obtain that 
\[
\frac{2L_p}{\vartheta} < \nu_{k} < \frac{2L_p}{\vartheta}, 
\]
which is impossible. Hence no index $k \in
\iibe{k_0}{N_1(\epsilon)-1}$ exists such that $\nu_k
<2L_p/\vartheta$. Thus, $k_1 \leq k_0$ exists by definition of $k_1$
in \req{k1def}. By the same definition, we finally deduce that $\nu_k
\geq 2L_p/\vartheta$ for all $k\geq k_1$, in turn implying
\req{sigmabig} because of \req{sigkupdate}. 
}

\noindent
Observe that \req{kstardef} depends on the ratio
$L_p/\nu_0$ which is the fraction by which $\nu_0$
underestimates the Lipschitz constant. This lemma thus implies that
the probability of $\calK_1$ is at least the probability that $k_*$
$\omega$-true iterations are performed, which we now investigate under
the following assumption.

\vskip 5pt
\textbf{AS.6} 
There exists an $\omega\in (0,1)$ and a $\pione>0$ such that for $S_k$
drawn randomly, 
$$
\prob\Big[ \mbox{iteration $k$ is $\omega$-true} \, \mid \, x_k=\bar x_k, \, \sigma_k=\bar \sigma_k \Big]\ge \pione,
$$
for any $\bar x_k\in \Re^n$, any $\bar \sigma_k\in [\vartheta \nu_0, \sigma_{\max}]$ and any $k\in \iiz{N_1(\epsilon)-2}$, where $\prob[X]$ denotes the probability of the event $X$.
Moreover, the occurrence of $k$-th iteration being $\omega$-true is conditionally independent of 
the occurrence of iterations $0, \ldots, k-1$ being $\omega$-true given $x_k=\bar x_k$ and $\sigma_k=\bar \sigma_k$.
\vskip 5pt

\noindent

This assumption differs from Assumption~1 in \cite[page~71]{Shao21} in that 
it now it makes the probability of an $\omega$-true iteration conditional 
not only on $x_k$ but also on $\sigma_k$, which we feel is reasonable given 
the isotropic nature of the regularisation term in \req{model}.
Note that a suitable value for $\omega$ may depend on the bounds on 
$\sigma_k$ (as we will see below in Lemmas~\ref{useful1}, \ref{useful2} and 
\ref{useful1B}). Assumption AS.6 can be ensured by 
suitably using Johnson-Lindenstrauss embeddings \cite{DasgGupt02} and results are 
available in \cite{Shao21}  for $p\in\{1,2\}$. We will analyse such cases 
in Section \ref{select-s}.

Before using AS.6 and $\pione$ directly, we first recall a known
probabilistic result.

\vskip 5pt
\llem{LemmaShao}{
For all nonnegative $i$, let $\calA_i$ be an event which can be true or
false and is conditionally independent of $\calA_0, \calA_1, \ldots
\calA_{i-1}$. For any $\bar x_i\in \Re^n$ and $\bar \sigma_i\in
     [\vartheta \nu_0, \sigma_{\max}]$, suppose that
     $\prob\Big[\calA_{i} \mbox{ is true} \, \mid \, x_i=\bar x_i, \,
       \sigma_i=\bar \sigma_i \Big]\ge \pi$, with $\pi\in (0,1)$. For
     $k\geq 0$, let $\calW_k=\{i \in \iiz{k-1} \mid \calA_i \mbox{ is
       true}\}$. 
Then, for any given $\delta_1\in (0,1)$,
\beqn{NTShao}
\prob\Big[|\calW_k| > (1-\delta_1)\pi k\Big] \ge 1-e^{-\frac{\delta_1^2}{2}\pi k}.
\eeqn
}
\proof{See \cite[Lemma 4.3.1]{Shao21} where, as mentioned above,
we now consider the ``state'' of the algorithm at iteration $i$ to 
comprise both $x_i$ and $\sigma_i$.}

\noindent
We are now in position to use this result to obtain a lower bound on
the probability that $k_*$ $\omega$-true iterations are performed, and 
that $k_1$ is well-defined.

\llem{Nits}{
Suppose that AS.1, AS.3 and AS.6 hold and let $\delta_1 \in (0,1)$ be
given. Let
\beqn{kdiamonddef}
k_\diamond \eqdef \left\lceil\frac{k_*}{(1-\delta_1)\pione}\right\rceil,
\eeqn
where $k_*$ is given by \req{kstardef}. Then 
\beqn{NTitsn}
\prob\Big[ \calK_1 \mid N_1(\epsilon) > k_\diamond \Big]
\ge  1- e^{-\frac{\delta_1^2}{2}\pione k_\diamond} \eqdef \pioneone.
\eeqn
}

\proof{
Identifying $\calA_i = \{\mbox{iteration $i$ is } \omega\mbox{-true}\}$,
  Lemma \ref{LemmaShao} with $\pi = \pione$ and $k_0=k_\diamond$ gives
  that the probability that at least $k_*$ $\omega$-true iterations have been
  performed during iterations $0$ to $k_\diamond-1$ is at least
  $\pioneone$.  The desired conclusion then
  follows from  Lemma~\ref{sigmanotbig}.
}

\noindent
We finally propose a variant of the well-known ``telescoping sum''
argument adapted to our probabilistic setting to 
derive the desired evaluation complexity bound.

\lthm{complexity}{
Suppose that AS.1, AS.2, AS.3, AS.4, AS.5 and AS.6 hold, 
that $\delta_1\in (0,1)$ is given and that the \skoffarp\ algorithm
is applied to problem \req{problem}. Define
\beqn{kapoff1}
\kap{SKOFFARp}
\eqdef \frac{4\left[L_p +(p+1)!(f_{\max}-f_{\rm low}) \right]}
            {\vartheta\nu_0\omega^{\sfrac{p+1}{p}}(1-\delta_1)\pione}
\eeqn
where $f_{\max}$ is defined in \req{fmaxdef}.
Then
\beqn{complexity1}
\prob\Big[ N_1(\epsilon) \leq \kap{SKOFFARp} \, \epsilon^{-\sfrac{p+1}{p}}+4 \Big]
\geq \left(1 - e^{-\frac{\delta_1^2}{2}\pione k_\diamond}\right)^2
\eeqn
where $k_\diamond$ is defined by \req{kdiamonddef}.
}

\proof{
  First note that \req{kstardef} and \req{kdiamonddef} imply that
  \beqn{kdiamondbound}
  k_\diamond \leq\frac{1}{(1-\delta_1)\pione}
  \left(\frac{2L_p}{\vartheta\nu_0\omega^{\sfrac{p+1}{p}}}\right)\epsilon^{-\sfrac{p+1}{p}}+1.
  \eeqn
Thus, given \req{kapoff1},
\beqn{P1}
\prob\Big[N_1(\epsilon)
\leq \kap{SKOFFARp} \, \epsilon^{-\sfrac{p+1}{p}}+4\, \mid \, N_1(\epsilon) \leq 2 k_\diamond + 2 \Big]
= 1.
\eeqn

Suppose now that $N_1(\epsilon) > k_\diamond + 2 > k_\diamond$
and that $\calK_1$ occurs. Consider an iteration $j >
k_\diamond\geq k_1$ (note that $k_1$ is well-defined) such that $j+1 < N_1(\epsilon)$
and suppose furthermore that iteration $j$ is $\omega$-true, a
situation which occurs with probability at least $\pione$ because of AS.6.
From the fact that $\calK_1$ occurs, $N_1(\epsilon)>k_\diamond$ and 
the definition of $k_1$ in \req{k1def}, we have that $\sigma_j \geq 2L_p$ 
and we may apply Lemma~\ref{Ifsigmabig}, yielding \req{sigma-upper} 
for iteration $j$. Since this iteration is also $\omega$-true, \req{sigma-upper} 
and inequality \req{crucial} also hold for iteration $j$.
Moreover, the fact $\calK_1$ occurs ensures (because of
Lemma~\ref{uppersigmak}, \req{sigkupdate}, the non-decreasing
nature of $\nu_k$ and the identity $\sigma_0=\nu_0$) that
$\sigma_j\in [\vartheta \sigma_0,\sigma_{\max}]$.
Finally, $\|g_{j+1}\|\geq \epsilon$ because $j+1 < N_1(\epsilon)$.
Combining these observations, we obtain that
\beqn{inproof1}
f(x_j) - f(x_{j+1})
\geq \frac{\sigma_j \|s_j \|^{p+1}}{2(p+1)!}
\geq \frac{\sigma_j \omega^\sfrac{p+1}{p} \|g_{j+1}\|^\sfrac{p+1}{p}}{2(p+1)!} 
\geq\frac{\vartheta\nu_0\omega^\sfrac{p+1}{p} \epsilon^\sfrac{p+1}{p}}
          {2(p+1)!}  
\eeqn
with probability (conditional to $\calK_1$ and
$N_1(\epsilon)>k_\diamond+2$) at least $\pione$.
Applying now Lemma~\ref{LemmaShao} to iterations of index $k_\diamond+1$ to $j$ with
\[
\calA_{i-k_\diamond} = \{\tim{\req{inproof1} holds at iteration $i-k_\diamond$} \},
\ms
\pi = \pione \tim{and} k = j-k_\diamond,
\]
we deduce that, for all $j\in\iibe{k_\diamond+1}{N_1(\epsilon)-2}$,
\[
\prob\left[ |\calV_j| \geq (j-k_\diamond)(1-\delta_1)\pione \, | \, \calK_1 \tim{and} N_1(\epsilon)>k_\diamond\right]
\geq 1 - e^{-\frac{\delta_1^2}{2}\pione (j-k_\diamond)}
\]
where $\calV_j \eqdef \{ i \in \iibe{k_\diamond+1}{j} \mid
\mbox{\req{inproof1} holds at iteration $i$} \}$. 
In particular, we have that
\beqn{inproof2}
\prob\left[ |\calV_j| \geq (j-k_\diamond)(1-\delta_1)\pione \, | \, \calK_1 \tim{and} N_1(\epsilon)>2k_\diamond+2\right]
\geq \pioneone,
\eeqn
with $\pioneone$ defined in \req{NTitsn}, for all $j\in \iibe{2k_\diamond+1}{N_1(\epsilon)-2}$.
We also know from
Lemma~\ref{Ifsigmabig} and the definition of $k_1$ in \req{k1def} that
the sequence $\{f(x_j)\}$ is non-increasing for $j \geq k_1$, and thus
that 
\[
f(x_{k_1}) - f(x_{j+1})
= \bigsum_{i=k_1}^{j}[f(x_i) - f(x_{i+1})]
\geq \bigsum_{i=k_\diamond+1}^{j}[f(x_i) - f(x_{i+1})]
\geq |\calV_j| \min_{i\in \calV_j } [f(x_i) - f(x_{i+1})].
\]
Combining this inequality with \req{inproof1} and \req{inproof2} then
yields that
\[
\prob\left[ 
f(x_{k_1}) - f(x_{j+1}) \geq (j-k_\diamond)(1-\delta_1)\pione\;\kappa_2^{-1}\epsilon^\sfrac{p+1}{p}
\, \mid \, \calK_1 \tim{and} N_1(\epsilon)>2k_\diamond +2\right] 
\geq \pioneone
\]
where
\beqn{kap2def}
\kappa_2 = \frac{2(p+1)!}{\vartheta\nu_0\omega^\sfrac{p+1}{p}},
\eeqn
and thus, because of AS.3, that
\[
\prob\left[ 
f(x_{k_1}) - f_{\rm low} \geq \kappa_2^{-1}(1-\delta_1)\pione (j-k_\diamond)\,\epsilon^\sfrac{p+1}{p}
\, \mid \, \calK_1 \tim{and} N_1(\epsilon)>2k_\diamond +2\right]
\geq \pioneone.
\]
Furthermore, \req{fmaxdef} in  Lemma~\ref{fk1boundsigmak1} then implies that
\[
\prob\left[
j-k_\diamond \leq \frac{\kappa_2}{(1-\delta_1)\pione}(f_{\max} - f_{\rm low})\,\epsilon^{-\sfrac{p+1}{p}}
\, \mid \, \calK_1 \tim{and} N_1(\epsilon)>2k_\diamond +2\right]
\geq  \pioneone,
\]
Since $j$ is arbitrary between $2k_\diamond+1$ and $N_1(\epsilon)-2$, we obtain that
\[
\prob\left[
  N_1(\epsilon)\leq \bigfrac{\kappa_2}{(1-\delta_1)\pione}
  (f_{\max} - f_{\rm low}) \,\epsilon^{-\sfrac{p+1}{p}}+ k_\diamond +2
\, \mid \, \calK_1 \tim{and} N_1(\epsilon)>2k_\diamond+2\right]
\geq \pioneone,
\]
which, given the definitions of $\kappa_2$ in \req{kap2def},  of
$\kap{SKOFFARp}$ in \req{kapoff1} and inequality \req{kdiamondbound},
yields that 
\[
\prob\left[
  N_1(\epsilon)\leq \kap{SKOFFARp} \, \epsilon^{-\sfrac{p+1}{p}}+4
\, \mid \, \calK_1 \tim{and} N_1(\epsilon)>2k_\diamond+2\right]
\geq \pioneone.
\]
Therefore, from \req{P1},  the fact that
\[
\prob\Big[\calK_1 \mid N_1(\epsilon)>2k_\diamond+2 \Big]
\ge \prob\Big[\calK_1 \mid N_1(\epsilon)>k_\diamond \Big]
\]
and Lemma~\ref{Nits}, we finally obtain that
\[
\begin{array}{l}
\prob\left[
  N_1(\epsilon)\leq \kap{SKOFFARp} \, \epsilon^{-\sfrac{p+1}{p}}+4 \right] \\*[2ex]
\hspace*{1cm}
= \prob\left[
  N_1(\epsilon)\leq \kap{SKOFFARp} \, \epsilon^{-\sfrac{p+1}{p}}+4
 \,  \mid \, \calK_1 \tim{and} N_1(\epsilon)>2k_\diamond+2\right] \\*[2ex]
\hspace*{2cm}
\times \,
\prob\Big[\calK_1 \mid N_1(\epsilon)>2k_\diamond +2\Big] \,
\times \, \prob\Big[N_1(\epsilon)>2k_\diamond +2\Big]
+ \, 1 \,\times \, \prob\Big[N_1(\epsilon)\leq 2k_\diamond +2\Big]\\*[2ex]
\hspace*{1cm}
\geq \prob\left[
  N_1(\epsilon)\leq \kap{SKOFFARp} \, \epsilon^{-\sfrac{p+1}{p}}\!+4
 \,  \mid \, \calK_1 \tim{and} N_1(\epsilon)>2k_\diamond+2\right] 
\times \,
\prob\Big[\calK_1 \mid N_1(\epsilon)>k_\diamond\Big]  \\*[2ex]
\hspace*{1cm}
\geq (\pioneone)^2.
\end{array}
\]
Substituting the values of  $\pioneone$ given by \req{NTitsn} in this 
inequality then yields \req{complexity1}.
} 

\noindent
We now comment on this result.
\begin{enumerate}
\item As in the methods of \cite{Shao21} and  \cite{CartFowkShao22}, it is 
not necessary to evaluate the full-space derivatives 
$\{\nabla_x^j f(x_k)\}_{j=1}^p$ because only their sketched versions 
$\{\nabla_x^j f(x_k)[S_k \cdot]^j\}_{j=1}^p$ are used. As a consequence, 
the cost of evaluating the derivatives (not to mention that of computing the 
step) is potentially reduced  typically by a significant factor $\ell/n$.
We discuss below whether this advantage may be offset by the 
choice of $\omega$ in AS.6.

\item Because it is proved in \cite[Theorem~3.12]{GratJeraToin22c}
  that the $\calO(\epsilon^{-(p+1)/p})$ order bound for finding
  $\epsilon$-approximate critical points is sharp for the
  \al{OFFAR$p$} algorithm, the same is also true for
  Theorem~\ref{complexity} above, because \skoffarp\ 
  subsumes\footnote{The different conditions on the regularisation 
  parameter $\sigma_k$ only result in differences in the constants.} 
  \al{OFFAR$p$} if $S_k = I$ for all $k$.
  \item Considered as a worst-case evaluation complexity bound for $p=2$, 
  the order bound $\calO(\epsilon^{-3/2})$ is known to be optimal for a 
  large class of methods using first- and second-derivatives 
  \cite{CartGoulToin18a}, justifying the title of this paper.
\item Note that \req{kstardef} and \req{kdiamonddef} not only
  imply \req{kdiamondbound}, but also that $k_\diamond$ is at least a
  (significant) fraction of $\epsilon^{-(p+1)/p}$, which, for
  meaningul values of $\epsilon$, is a reasonably large number.  Moreover,
  $(k_\epsilon-k_\diamond)$ is expected to be at least of the same order.
  Thus the factor
  \[
  \left(1 - e^{-\frac{\delta_1^2}{2}\pione k_\diamond}\right)
  \]
  in the right-hand side of \req{complexity1} is expected to be
  very close to 1.
\item The parameter $\delta_1$, which we are still free to choose in
  (0,1) occurs in \req{kapoff1} and in the exponentials of
  \req{complexity1}. A quick calculation indicates that choosing
  $\delta_1$ close to 1 can improves the bound on the
  right-hand side of \req{complexity1} (although marginally because of
  our previous comment) while its possibly detrimental
  effect on \req{kapoff1} occurs because of the factor
  $1/(1-\delta_1)$ which must be kept bounded. Given the magnitude of
  the other factors in these formulae, values such as
  $\delta_1=\half$ or $\delta_1=\sfrac{1}{10}$ could be
  considered acceptable.
\item As can be expected, the conditions for a random embedding given
  by \req{true-1} and AS.6 have a significant impact on
  the result, which significantly degrades if $\omega$ and/or
  $\pione$ tends to zero. 
\item As we have mentioned above, the objective function is not evaluated by the
  \skoffarp\ algorithm and  the trial point $x_k+s_k$ is
  always accepted as the next iterate. Thus no
  distinction is necessary in the stochastic analysis between
  ``successful'' iterations (where the step is accepted because the
  objective function has decreased enough) and  ``unsuccessful'' ones.
  This distinction had however to be taken into account in the
  analysis of \cite{Shao21} for more standard trust-region and
  adaptive-regularisation methods using functions values, leading to several
  different types of iterations whose numbers have to be bounded.   
\end{enumerate}

\section{Selecting random subspaces}\label{select-s}

We now turn to ways in which $\omega$-true iterations can be shown to
happen with suitable probability $\pione$, thereby satisfying AS.6. A
natural approach is to rely on Johnson-Lindenstrauss embeddings 
and results are available in the literature for 
$p\in\{1,2\}$. Restricting ourselves to such values of $p$  and using
\cite[Definition 5.3.1]{Shao21} (see also \cite{Wood14}, for instance), we say
that, for some given ``preservation parameter'' $\alpha_S \in (0,1)$
and for some positive scalar $S_{\max}$ independent of $k$,
iteration $k$ is $(\alpha_S,S_{\max})$-embedded whenever
\beqn{Smax-def}
\|S_k\|\le S_{\max},
\eeqn
and for
\beqn{Mk-def}
M_k \eqdef [g_k, H_k] \in \Re^{n\times n+1},
\eeqn
we have that
\beqn{true-1}
\|S_k M_k z\|
\geq \alpha_S \|M_k z\|
\tim{ for all } z \in \Re^{n+1},
\eeqn
where $H_k=\nabla^2_s f(x_k)$ if $p=2$ and $H_k=0_{n\times n}$ if $p=1$.
This condition is said to define a \textit{one-sided random embedding}
of the second-order Taylor's series.

Given such a one-sided random embedding, we now adapt an argument of \cite{GratToin21}
and verify that \req{crucial} holds at $(\alpha_S,S_{\max})$-embedded iterations.

\llem{useful1}{
  Suppose that $p\in \{1,2\}$, $p!\le L_p$, that AS.1, AS.2, AS.3, AS.4 and AS.5 hold and 
  that iteration $k\geq 0$ of the \skoffarp\ algorithm is 
  $(\alpha_S,S_{\max})$-embedded (in the sense of \req{true-1}). Then
  \beqn{iscrucial}
  \|s_k\|^p \geq \frac{p!\,\alpha_S}{\alpha_S L_p + \theta S_{\max}\sigma_{\max} } \| g_{k+1}\|.
  \eeqn
  Thus iteration $k\in\iiz{N_1(\epsilon)-2}$ is $\omega$-true (in the sense of \req{crucial}) 
  with $\omega = \frac{p!\,\alpha_S}{\alpha_S L_p + \theta S_{\max}\sigma_{\max}}$.
}

\proof{
First note that applying the chain rule gives that
\[
\nabla^1_{\widehat s}  \widehat T_{f,p}(x_k,\widehat s_k)
= S_k \nabla^1_s  T_{f,p}(x_k,s_k)
= S_k(g_k + H_k s_k) = S_k M_k (1, s_k^T)^T
\]
and, since the iteration $k$ is $(\alpha_S,S_{\max})$-embedded, \req{true-1} gives
that
\[
\|\nabla^1_{\widehat s}  \widehat T_{f,p}(x_k,\widehat s_k)\|
\geq \alpha_S \| M_k (1, s_k^T)^T\|
= \alpha_S \|\nabla^1_s  T_{f,p}(x_k,s_k)\|.
\]
Condition \req{gradstep}, the definition $s_k=S_k^T\widehat s_k$ and \req{Smax-def} then yield that
\beqn{bound_nabla_T}
\|\nabla_s^1 T_{f,p}(x_k,s_k)\|
\le\frac{\|\nabla^1_{\widehat s}\widehat T_{f,p}(x_k,\widehat s_k)\|}{\alpha_S}
\le\frac{\theta \frac{\sigma_k}{p!}S_{\max}\|S_k^T\widehat s_k\|^p}{\alpha_S}
\le \frac{\theta S_{\max}\sigma_k}{p!\,\alpha_S}\|s_k\|^p.
\eeqn
Successively using the triangle inequality, condition
\req{bound_nabla_T} and \req{Lip-g} (for $p\in\{1,2\}$), we deduce that
\[
\|g_{k+1}\| \leq \|g_{k+1}- \nabla_s^1 T_{f,p}(x_k,s_k)\| + \|\nabla_s^1 T_{f,p}(x_k,s_k)\|
\leq \bigfrac{1}{p!} L_p \|s_k\|^p +\frac{\theta S_{\max}\sigma_k}{p!\,\alpha_S}\|s_k\|^p.
\]
The inequality \req{iscrucial} follows by rearranging the terms and
using the bound \req{sigmamaxdef} in Lemma~\ref{uppersigmak}.
That iteration $k$ is $\omega$-true for $k\in\iiz{N_1(\epsilon)-2}$ follows from the fact that, by definition,
$\|g_{k+1}\|\geq \epsilon$ for these values of $k$.
}

\noindent
Although this lemma essentially recovers the result of
\cite[Lemma~5.3.2]{Shao21}, its proof is considerably simpler.
Note that \req{iscrucial} is significantly stronger than \req{crucial}, suggesting that
\req{true-1} might itself be stronger than necessary.
Also observe that we could replace condition \req{gradstep} by the more permissive
\[
\|\nabla^1_{\widehat s} \widehat T_{f,p}(x_k,\widehat s_k)\|
\leq \theta \frac{\sigma_k}{p!} \|S_k\|\,\|S_k^T\widehat s_k\|^p
\]
or
\[
\|\nabla^1_{\widehat s} \widehat T_{f,p}(x_k,\widehat s_k)\|
\leq \theta \frac{\sigma_k}{p!} \kappa_{S,k}\,\|S_k^T\widehat s_k\|^p
\]
without altering the above theory, but at the price of computing $\|S_k\|$ 
or estimating a uniform bound on $\kappa_{S,k}$ (such as $S_{\max}$).

It is also possible to apply Shao's approach to
``sparse Hessians'' (for $p=2$) as follows. For some constants $(\alpha_S,\gamma_S)$ 
such that $\alpha_S\in(0,1)$ and $\gamma_S\in [0,2\alpha_S)$ and $S_{\max}>0$, we now (re)define iteration $k$ to be
$(\alpha_S,\gamma_S, S_{\max})$-embedded whenever
\beqn{trues-1}
\|S_k\|\le S_{\max}, \quad 
\|S_k g_k \| \geq \alpha_S \|g_k\|
\tim{ and }
\|S_k H_k\| \leq \sqrt{\gamma_S \|g_{k+1}\|}.
\eeqn

We then obtain the following result based on \cite[Lemma~5.4.1]{Shao21}.

\llem{useful2}{
  Suppose that AS.1 and AS.3 hold and that, for a particular
  realization, iteration $k\geq 0$ of the \skoffar\al{2} algorithm
  is $(\alpha_S,\gamma_S, S_{\max})$-embedded (in the sense of \req{trues-1}). 
  Then \req{crucial} holds and iteration $k$ is $\omega$-true. 
}

\proof{  
Let $a = \|S_kH_k\|$. Then \req{gradstep} gives that
\[
\alpha_S\|g_k\|
\leq \|S_kg_k\|
\leq \|S_k(g_k+H_ks_k)\| + \|S_kH_ks_k\|
\leq  \half \theta S_{\max} \sigma_k \|s_k\|^2 + a\|s_k\|,
\]
and therefore, using the triangle inequality,
\req{Lip-g} and the fact that iteration $k$ is $(\alpha_S,\gamma_S, S_{\max})$-embedded,
\[
\alpha_S \|g_{k+1}\|
\leq \alpha_S \|g_{k+1}-g_k\| + \alpha_S \|g_k\|
\leq \half \alpha_S L_2 \|s_k\|^2 +\half \theta S_{\max} \sigma_k \|s_k\|^2 + a \|s_k\|.
\]
Defining $b=\alpha_S L_2 + \theta S_{\max} \sigma_k$, we obtain that
\[
\|s_k\|^2 + \left(\frac{2a}{b}\right)\|s_k\|-\frac{2\alpha_S\|g_{k+1}\|}{b} \geq 0,
\]
yielding that
\[
\left(\|s_k\| + \frac{a}{b}\right)^2
\geq \frac{2\alpha_S\|g_{k+1}\|}{b} + \left(\frac{a}{b}\right)^2
\]
and thus that
\[
\|s\| \ge \sqrt{\frac{2\alpha_S\|g_{k+1}\|}{b} + \left(\frac{a}{b}\right)^2} - \frac{a}{b}.
\]
Assuming, without loss of generality, that $b=\alpha_S L_2 + \theta S_{\max} \sigma_k\ge 1$, we deduce that
\[
\|s\|
\ge \frac{1}{b}\left[\sqrt{2\alpha_S\|g_{k+1}\| + a^2} - a\right]
\]
Since the function $\sqrt{c+ t^2}-t$ (for $c>0$) is decreasing as a
function of $t\ge0$ and since
$a = \|S_kH_k\| \leq \sqrt{\gamma_S \|g_{k+1}\|}$ because iteration $k$ is $(\alpha_S,\gamma_S)$-true, we deduce
that
\[
\begin{array}{lcl}
\|s\|
& \ge &\bigfrac{1}{b}\left[\sqrt{2\alpha_S\|g_{k+1}\| + \gamma_S\|g_{k+1}\|}
  - \sqrt{\gamma_S\|g_{k+1}\|}\right]\\*[2ex]
& \ge & \bigfrac{\sqrt{2\alpha_S}-\sqrt{\gamma_S}}{\alpha_S L_2 + \theta S_{\max} \sigma_k}\sqrt{\|g_{k+1}\|}\\*[2ex]
& \ge & \bigfrac{\sqrt{2\alpha_S}-\sqrt{\gamma_S}}{\alpha_S L_2
  + \theta S_{\max}}\sigma_{\max} \sqrt{\|g_{k+1}\|},
\end{array}
\]
where we again used Lemma~\ref{uppersigmak} to obtain the last inequality.
} 

\noindent
Thus, an $(\alpha_S,\gamma_S, S_{\max})$-embedded iteration (in the sense of
\req{trues-1}) is $\omega$-true (in the sense of \req{crucial}) for
$\omega =( \sqrt{2\alpha_S}-\sqrt{\gamma_S})/(\alpha_S L_2
+ \theta S_{\max}\sigma_{\max})$.  
Also notice that, should we replace \req{trues-1} by
\beqn{trues2-1}
\|S_k\|\le S_{\max}, \quad
\|S_k g_k \| \geq \alpha_S \|g_k\|
\tim{ and }
\|S_k H_k\| \leq \sqrt{\gamma_S \epsilon}
\tim{ for }
k < N_1(\epsilon)-1,
\eeqn
then the definition of an $(\alpha_S,\gamma_S, S_{\max})$-embedded iteration is closer
to that of \cite{Shao21}, obviously ensuring \req{trues-1} with a
right-hand side of  its third part now independent of $S_k$.

The reader may now recall that AS.6 states that \req{crucial},
\req{true-1}, \req{trues-1} or \req{trues2-1} (or the first part of
\req{trues-1} or \req{trues2-1}) should hold at iteration $k$ with
positive probability $\pione$. In
\cite[Lemma~5.3.1]{Shao21} or \cite[Theorem~2.3]{Wood14} (see also
\cite[Lemma~3.1]{VuPoirDambLibe19}) it is argued that by choosing $\calS$ to be
the distribution of $\ell \times n$ scaled Gaussian matrices, \req{true-1} holds with
probability
\beqn{bad-pione}
\pione = 1 - e^{-\frac{\ell(1-\alpha_S)}{C_\ell} + {\rm rank}(M_k)},
\eeqn
where $C_\ell> \sfrac{1}{4}$ is an absolute constant. 

Unfortunately, the expression \req{bad-pione}  requires that
\beqn{low-rank}
{\rm rank}(M_k) < \frac{\ell(1-\alpha_S)}{C_\ell}, 
\eeqn
thereby limiting the applicability of the result for $p>1$ when considering general problems with full-rank Hessians. But this can be acceptable for a class of problems with low-rank Hessian, as we illustrate in Section~\ref{numerics-s}.
Satisfying the third part
of \req{trues2-1} with positive probability is also possible when $H_k$ is
very sparse, which also imposes a significant restriction.
Other choices for the distribution exist, such as hashing,
scaled hashing, sampling matrices, or ``fast Lindenstrauss transforms''
(see \cite[Chapter~2]{Shao21} or \cite[page~16]{Wood14}). Although
possibly more economical in terms of algebraic operations, they
appear to suffer from the same geometric precondition: their number of rows
$\ell$ should be of the order of the Hessian's rank, which is
problematic for the general case where the Hessian is full rank. 
However, note that ${\rm rank}(M_k) = 1$ when $p=1$, essentially
avoiding this problem, making the first-order variant of the
algorithm applicable to a much larger class of problems.

Should one be ready to trade the optimal complexity for getting rid
of the low-rank requirement, an algorithm
using quadratically regularised quadratic models with inexact Hessians
can also be defined and analysed (see Appendix). Under suitably
modified assumptions, the evaluation complexity of this algorithm can
be shown to be of order $\calO(\epsilon^{-2})$, matching the
theoretical results of \cite{CartFowkShao22} for a random subspace
version of the adaptive-regularisation algorithm \textit{using function
values}. Unfortunately, our numerical experience matches the
cautious conclusions of this reference, which is why we do not
investigate it further.

Finally, note that the constant \req{kapoff1} involves $S_{\max}^{\frac{p+1}{p}}$ 
due to its dependence on $\omega^{\frac{p+1}{p}}$.
In the case of scaled Gaussian matrices, we know that
\beqn{Gaussian-Smax}
S_{\max} \leq  \beta \eqdef 1.5 + \sqrt{n/\ell}
\eeqn
with high probability for the values considered of $\delta_1$ (see
\cite[Lemma~4.4.4]{Shao21} for instance), resulting in a dependence  
of the constant \req{kapoff1} on $(n/\ell)^{\frac{p+1}{2p}}$.  For $p
= 1$, this offsets, complexity-wise, the benefit of cheaper gradient
evaluations by a factor of $\ell/n$, but the complexity bound is
rarely tight and savings in gradient evaluations are sometimes
possible in practice. For $p=2$, the advantage of cheaper gradients
(assuming \req{low-rank}) increases compared to $p=1$ because the
  denominator of \req{kapoff1} now depends on $(n/\ell)^{3/4}$. We
  show in Section~\ref{numerics-s} that this theoretical advantage
  translates into significantly better numerical behaviour. Moreover,
  the cost of computing the step $s_k$ is reduced by the decrease in
  dimension of the linear system. The advantage of cheaper gradients
grows when $p$ grows (and the method is applicable). 

\section{Numerical illustration}\label{numerics-s}

We now numerically illustrate the behaviour of
\skoffarduedot, the second-order version of \skoffarp. 
We report results obtained using Matlab R2024a for 14 problems from the {\sf CUTE}st test problems \cite{GoulOrbaToin15b} as provided in Matlab by OPM \cite{GratToin21b}.
All problems except {\tt arglina} and {\tt tridia} are nonconvex.
The original dimension of the problem, say $\hat n$, was enlarged using the affine transformation $x=A \hat x $, $\hat x\in R^{\hat n}$, $x\in R^n$, $n \ge \hat n$,  $A\in R^{n \times \hat n}$ being an orthonormal matrix generated by the (Matlab supplied) discrete cosine transform, therefore yielding problems with Hessians of rank at most $\hat n$ and ensuring \req{low-rank} when $n$ grows. 

We ran a Matlab implementation of a modified version 
of the \skoffarp~where we defined $\calS$ to be the distribution of 
$\ell \times n$ scaled Gaussian matrices.
The first modification is identical to that described in \cite{GratJeraToin22c} for the 
\al{OFFAR$p$} algorithm, in that \req{sigkupdate} is replaced by
\[
\sigma_k = \max[ \vartheta \nu_k, \xi_k \mu_k ]
\]
where $\xi_k \in (0,1)$ is an adaptive scaling parameter (see \cite{GratJeraToin22c} for details)
and where $\mu_k$ is defined by \req{muk-def} with $\mu_{-1}= 10^3$. The second change
avoids the (potentially very) costly computation of $\|S_k\|$ by using
$\kappa_{S,k}= \beta$ as given by \req{Gaussian-Smax}. This change was
made after running the more expensive code using $\kappa_{S,k}=\|S_k\|$ as suggested by
\req{Lpmuk} on a few problems and observing that the results obtained with the theoretically 
weaker $\kappa_{S,k} = \beta$ did not degrade the code efficiency, if at all.
The regularised quadratic was minimized approximately ($\theta = 1.01(1 + \sqrt{n/\ell})$ 
using a Lanczos-based solver for such functions (see \cite[Section~10.2]{CartGoulToin22}). 
We also chose $\vartheta = 10^{-3}$ and terminated the optimization as soon as the threshold
$\|g_k\| \leq 10^{-3}$ was reached. 
All computations were performed on a Dell Precision laptop with 16 cores 
at 2.6 GHz and 62.5 GB of memory, running Matlab 2024a with Ubuntu 20.04.6 LTS.  

For each run, we computed the number of gradient evaluations weighted to 
reflect the reduced evaluation cost in the subspace of dimension $\ell$. 
We counted the cost of evaluating a Hessian as the product of the dimension 
times the cost of evaluating one gradient (as happens for finite-difference 
approximations), the weighted cost of an iteration (now involving the computation of 
one gradient and one Hessian) then becoming $w_2(\tau,n) = (\tau+n\tau^2)/(1+n)$, 
where $\tau=\ell/n$. Thus this $w_2$ weighting reflects the cost of running 
the second-order \skoffardue with sketching parameter $\tau$ compared to 
running it in full-dimension. To take this into account, the 
maximum number of iterations was set to $10^5$ divided by this factor. 
Table~\ref{resp1} reports the average $w_2$-weighted iteration costs 
for \skoffardue to reach convergence for $n = 1000 \,{\hat n}$ and for 
decreasing values of the ratio $\tau$ from 1 (irrealistic) to $10^{-3}$, 
averaged over 10 independent runs.

\begin{table}[ht]{\small 
\begin{center}
\begin{tabular}{|l|r|rrrrrrr|}
\hline
        &                          & \multicolumn{7}{|c|}{\skoffardue}\\
Problem & \multicolumn{1}{|c|}{$n$} & $\tau=1$ & $10^{-1}$ & $5 \cdot 10^{-2}$ & $2.5\cdot 10^{-2}$& $10^{-2}$ &  $5.10^{-3}$ & $10^{-3}$ \\
\hline
{\tt arglina}   &  10000 & 19.6000 & 0.8358 & 0.3394 & 0.1399 & 0.0443 & 0.0188 & 0.0027 \\
\hline
{\tt arwhead}   &  10000 &  6.7000 & 0.1191 & 0.0363 & 0.0125 & 0.0033 & 0.0013 & 0.0002 \\
\hline
{\tt broyden3d} &  10000 &  7.0000 & 0.1441 & 0.0393 & 0.0161 & 0.0054 & 0.0020 & 0.0004 \\
\hline
{\tt chandheu}  &  10000 &  6.0000 & 0.1231 & 0.0436 & 0.0160 & 0.0062 & 0.0032 & 0.0006 \\
\hline
{\tt dixmaana}  &  12000 & 13.0000 & 0.4253 & 0.2429 & 0.1290 & 0.0510 & 0.0244 & 0.0041 \\
\hline
{\tt eg2}       &  10000 &  5.3000 & 0.1611 & 0.0614 & 0.0233 & 0.0068 & 0.0028 & 0.0004 \\
\hline
{\tt engval2}   &   3000 & 18.4000 & 0.4082 & 0.2493 & 0.1323 & 0.0558 & 0.0275 & 0.0055 \\
\hline
{\tt helix}     &  10000 & 30.8000 & 1.7726 & 0.8606 & 0.4182 & 0.1892 & 0.1234 & 0.0241 \\
\hline
{\tt kowosb}  &  10000 & 2715.9000 & 152.9072 & 64.4489 & 27.2542 & 8.8252 & 3.8179 & 0.6299 \\
\hline
{\tt nzf1}      &  13000 &104.9000 & 7.3581 & 3.2245 & 1.4262 & 0.4701 & 0.2006 & 0.0297 \\
\hline
{\tt rosenbr}  &  10000 & 98 .7000 & 7.4097 & 3.7827 & 1.6012 & 0.5338 & 0.2675 & 0.0474 \\
\hline
{\tt sensors}   &  10000 & 18.0000 & 0.8398 & 0.3399 & 0.1393 & 0.0443 & 0.0192 & 0.0029 \\
\hline
{\tt tridia}    & 10000 &  14.7000 & 0.5565 & 0.2367 & 0.0961 & 0.0272 & 0.0112 & 0.0029 \\
\hline
{\tt watson}   &  10000 &  44.0000 & 3.9886 & 1.9084 & 0.8106 & 0.2847 & 0.1342 & 0.0146 \\
\hline
\end{tabular}
\end{center}
\caption{\label{resp1}Using the \skoffardue algorithm: average $w_2$-weighted  number of iterations for varying ratio $\tau=\ell/n$.}
} 
\end{table}

The results in Table \ref{resp1} show that the use of random subspaces
can bring substantial benefits, as long as \req{low-rank} holds. In
this context, \skoffardue is reliable and efficient on nearly all
problems (except for {\tt kowosb}) and tested values of
$\tau$. Performance globally increases with decreasing values of
$\tau$; it appears to be best for the smallest value
$\tau=10^{-3}$. Limited to the set of problems considered, these
results show that sketching pays off handsomely in terms of gradient
evaluations when used with low-rank Hessians and second-order
models. This, admittedly, ignores the cost of linear algebra, which
increases because products with $S_k$ have to be computed but also
decreases because the calculation of the step in the subspace is
significantly cheaper than in the full space.  

The reader might wonder at this point how (sketched) second-order methods 
might compare to standard first-order algorithms, where the Hessian is not evaluated.  We 
attempt to clarify this question by comparing, in Table~\ref{resp3}, our results for 
\skoffardue with those obtained by the well-known objective-function-free 
\al{ADAGRAD}-\al{N}orm algorithm \cite{DuchHazaSing11,WardWuBott19} and the norm-wise 
variant of \al{ADAM}\footnote{Using the momentum discounting
  factor $\beta = 0.9999$. It failed to converge on most problems with
  smaller values.} \cite{KingBa15}.  To make the comparison fair, we
have re-weighted the iteration counts in order to make them relative
to a single gradient evaluation (as is case for one iteration of \al{ADAGRAD} and \al{ADAM}) 
by using $w_1(\tau,n) = \tau + n\tau^2$.  In this table, the string 
``$>$100000'' indicates that the maximum number of iterations was reached. 

\begin{table}[ht]{\small 
\begin{center}
\begin{tabular}{|l|r|r|r|rrrrrr|}
\hline
      &   &\multicolumn{1}{|c}{} &\multicolumn{1}{|c}{}   &\multicolumn{6}{|c|}{\skoffardue}\\
Problem & \multicolumn{1}{|c|}{$n$} & \al{ADAM-N}  & \al{ADAG-N}  & $10^{-1}$ & $5 \cdot 10^{-2}$ & $2.5\cdot 10^{-2}$& $10^{-2}$ &  $5.10^{-3}$ & $10^{-3}$ \\
\hline
{\tt arglina}   &  10000 & 125 & 126 & 8358 & 3394 & 1399 & 443 & 118 & 27 \\
\hline
{\tt arwhead}   &  10000 & 45 & 45 & 1191 &  363 &  125 &   33 & 13 & 2 \\
\hline
{\tt broyden3d} & 10000  & 40 & 40 &  1441 & 393 & 160 & 54 & 20 & 4 \\
\hline
{\tt chandheu}   &  10000 & 51 & 51 &  1231 & 435 & 160 & 62 & 32 & 6 \\
\hline
{\tt dixmaana}  &  12000 & 697 & 710 & 5104 & 2915 & 1549 & 612 & 293 & 47 \\
\hline
{\tt eg2}       &  10000 & 106 & 104 &  1611 & 614 & 233 & 68 & 28 & 4 \\
\hline
{\tt engval2}  &  3000 & $>100000$ & 19266 &  1225 & 748 & 397 & 168 & 83 & 16 \\
\hline
{\tt helix}   &  10000 & 26142 & 53907 & 17727 & 8607 & 4183 & 1892 & 1234 & 241 \\
\hline
{\tt kowosb}   &  10000 & 295 & 296 & 611781 & 257860 & 109043 & 35309 & 15275 & 2520 \\
\hline
{\tt nzf1}      &  13000 & 8335 & 10323 & 95662 & 41921 & 18540 & 6112 & 2608 & 387 \\
\hline
{\tt rosenbr}   &  10000 & 26748 & 56173 & 74104 & 37830 & 16013 & 5338 & 2675 & 474 \\
\hline
{\tt sensors}   &  10000 & 189 & 167 &  8393 & 3399 & 1393 & 443 & 192 & 29 \\
\hline
{\tt tridia}    & 10000 & 50 & 50 & 5566 & 2367 & 962 & 272 & 112 & 29 \\
\hline
{\tt watson}   &  10000 & $>100000$ & 15132 &  39889 & 19085 & 8107 & 2848 & 1342 & 146 \\
\hline
\end{tabular}
\end{center}
\caption{\label{resp3}Using the \skoffardue algorithm: average $w_1$-weighted  number of iterations for varying ratio $\tau=\ell/n$.}
}
\end{table}

\begin{figure}[!ht]
	\centering
	\includegraphics[width=10cm,height=8cm]{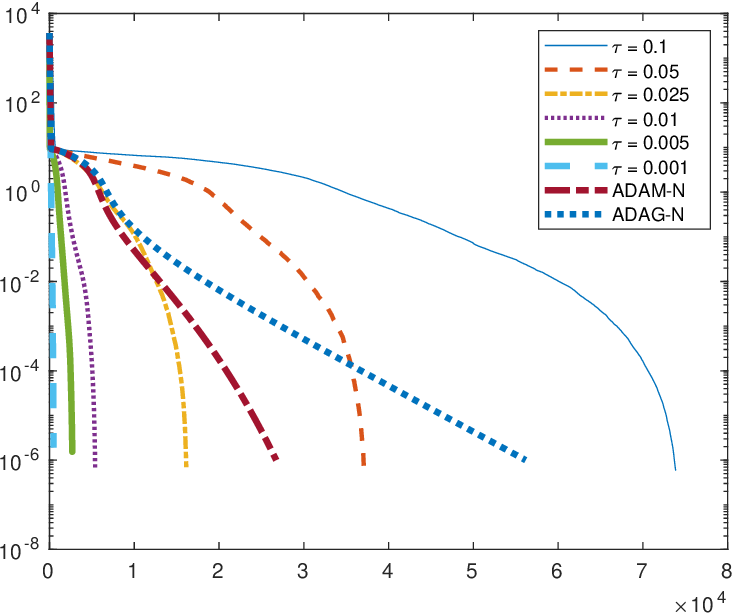}
	\caption{\label{fig}The behaviour of $f(x)$ when \al{ADAM-N}, \al{ADAG-N} and \skoffardue are run on {\tt rosenbr}, as a function of $w_1$-weighted iteration numbers, where \skoffardue uses $\tau =10^{-1}, 5 \cdot 10^{-2}, 2.5\cdot 10^{-2}, 10^{-2}, 5.10^{-3}$ and $10^{-3}$ (from right to left) }
\end{figure} 

These re-weighted results indicate that \textit{sketched second-order methods 
may be competitive with existing first-order algorithms when using a small 
$\tau$ for problems with low-rank Hessians}.  We also note the difference of 
performance between \al{ADAGRAD} and \al{ADAM}: while the latter may be more 
efficient when it works, it is less reliable than the former (as predicted 
by the theory). To provide further intuition, we also show, in 
Figure~\ref{fig}, how the value of the objective function evolves (although it 
is never computed in the course of the algorithm) for one instance of 
applying \al{ADAM}, \al{ADAGRAD} and \skoffardue to the {\tt rosenbr} problem, the latter with various choices of $\tau < 1$. Beyond the clearly 
faster convergence for smaller $\tau$, one also notices the concave nature of 
the curves, which contrasts with the convex curves one would often expect when 
using first-order methods. (The nearly vertical part of the curves between 
$10^4$ and $10^1$ corresponds to the first phase of minimization where all 
algorithms reach for the bottom of the steep, curving valley that is typical 
of {\tt rosenbr}.)  Because of this concavity, one also sees that requesting higher accuracy in \skoffardue is unlikely to require many more iterations.

\section{Conclusions and perspectives}\label{concl-s}

We have introduced an adaptive-regularisation algorithm for nonconvex unconstrained optimization that uses random subspaces and never computes the objective function's value, and have shown that its evaluation complexity is, in order, identical to that of the ``optimal'' adaptive full-space regularisation methods using function values.  The analysis covers finding approximate first-order critical points, but it is possible to extend the algorithm to ensure second-order criticality (along the lines of the MOFFAR algorithm in \cite{GratJeraToin22c}), albeit at the price of a very strong assumption on the recovery of the Hessian's minimum eigenvalue in random subspaces, a notoriously thorny problem (see \cite[Section~4.2.3]{Boma99}, for instance).  Our analysis also allows the use of models of arbitrary degree, but this generality may be of limited practical use since using degree higher than two appears to be mostly applicable to problems with low-rank, very sparse or partially separable derivatives.

Our theoretical and numerical results show that the approach is theoretically sound and that it can be significantly advantageous when its second-order variant is used
on problems with low-rank Hessians. 

{\footnotesize

\section*{Data availability}
  The test problems used in this study are available at 
  \url{https://github.com/gratton7/OPM}.

\section*{Acknowledgement}
  Ph. Toint  acknowledges the partial support of INDAM (through the GNCS grant CUP E53C22001930001) during his visit to the Department of the Industrial Engineering (DIEF) at the University of Florence in the fall 2023. He is also grateful for the continued and friendly support of DIEF and ANITI (Toulouse).

}

\appendix

\section{Quadratic regularisation for approximate second-order models}\label{GN-s}

We discuss here a context in which the low-rank assumption is
unnecessary and, motivated by \cite{CartFowkShao22}, consider using
quadratic regularisation in conjunction with approximate quadratic
models in which the Hessian  $\nabla_x^2 f(x)$ is approximated by a
positive-semidefinite symmetric matrix $B_k$.  At $x_k$, the
regularised model $\mkB(s)$ of $f(x_k+s)$ then takes the form 
\beqn{modelB}
\mkB(s) \eqdef \TB(x_k,s) +\frac{\sigma_k}{2} \|s\|^2,
\eeqn
with 
\beqn{taylor modelB}
\TB(x_k,s) \eqdef f(x_k)+\nabla f(x_k)^Ts+ \frac{1}{2} s^T B_k s.
\eeqn
To make the use of this model well-defined, we complete 
AS.2, AS.4 (for $p=1$) and AS.6  and make the following assumptions. 

\noindent
\textbf{AS.7} $f$ is continuously differentiable in $\Re^n$. 

\noindent
\textbf{AS.8} The gradient  of $f$ is globally Lipschitz
continuous, that is, there exist a non-negative constant $L_1$ such that
\[
\|\nabla_x^1 f(x) - \nabla_x^1 f(y)\| \leq L_1 \|x-y\| \, \text{ for all } x,y \in \Re^n.
\]

\noindent
\textbf{AS.9} The matrix $B_k$ is symmetric, positive-semidefinite and bounded for all $k \geq 0$, so that 
there exist a positive scalar $\kappa_B$ such that 
\beqn{kappaB}
\|B_k\|\le \kappa_B
\tim{ for } k\geq 0.
\eeqn
\vskip 5pt
Notice that AS.9\footnote{Alternatively, we could replace the condition that $B_k$ is positive-semidefinite by the weaker condition that $B_k+\sigma_k I$ is positive-semidefinite.} prevents the quadratic model \req{modelB} to be unbounded
below.  In particular, the use of the Gauss-Newton Hessian approximation for nonlinear least-squares problem is covered by AS.9, as well as the use of several quasi-Newton updating formulae.

Proceeding as in \skoffarp, we let $S_k$ be drawn from an iteration-independent distribution $\calS$ of $\ell \times n$ random matrices (with
$\ell<n$),  let $ s=S_k^T \widehat s$  be the full dimensional step and consider minimizing the 
sketched regularised model 
\beqn{modelB sub}
\widehat \mkB(\widehat s) \eqdef \widehat\TB(x_k,\widehat s) +\frac{1}{2}\sigma_k\|S_k^T \widehat s\|^{2},
\eeqn
where
\[
\widehat\TB(x_k,\widehat s) \eqdef f(x_k)+g_k^TS_k^T\widehat s+ \frac{1}{2} \widehat s^T S_k B_k S_k^T \widehat s.
\]
We note that, similarly to \req{equalmodels},
$\widehat \mkB(\widehat s)= \mkB( s)$.
The resulting \skoffar\al{2B} algorithm is stated  \vpageref{SKOFFAR2B}.

\algo{SKOFFAR2B}{OFFO adaptive regularisation with approximate second-order models (\skoffar\al{2B})}{
\begin{description}
\item[Step 0: Initialization: ] An initial point $x_0\in \Re^n$, a regularisation
parameter $\nu_0>0$ and a requested final gradient accuracy
$\epsilon \in (0,1]$ are given, as well as the parameters
$\theta > 1, \mu_{-1} \geq 0$  and  $0< \vartheta <1$. 
Set $k=0$.
\item[Step 1: Step calculation: ]
If $k=0$, set $\sigma_0=\nu_0$.
Otherwise, select a matrix $B_k$ satisfying AS.9 and
\[
\sigma_k  \in \Big[ \vartheta \nu_k, \max[\nu_k, \mu_k] \Big],
\]
where
\[
\mu_k = \max\left[ \mu_{k-1}, \frac{\|S_{k-1}g_k\| -
\|\nabla_s^1 \widehat \TB(x_k,s_k)\|}{\kappa_{S,k-1}.\|s_{k-1}\|}\right]
\]
with some $\kappa_{S,k-1}$ such that $\|S_{k-1}\| \leq \kappa_{S,k-1}$.
Draw a random matrix $S_k\in \Re^{\ell \times n}$ from
$\calS$ and compute a step $s_k=S_k^T \widehat s_k$  such
that $\widehat s_k$ sufficiently reduces  the
random model $\widehat \mkB$ defined in \req{modelB sub} in
the sense that  
\beqn{descentB}
 \widehat \mkB(\widehat s_k) - \widehat \mkB(0) <  0
\eeqn
and
\beqn{gradstepB}
\|\nabla_{\widehat s}^1 \widehat \TB(x_k,\widehat s_k)\|
 \leq  \theta \sigma_k\|S_k S_k^T\widehat s_k\|. 
\eeqn
\item[Step 2: Updates. ]
Set $x_{k+1} = x_k + s_k$ and
$\nu_{k+1} =  \nu_k + \nu_k \| s_k\|^2$.
Increment $k$ by one and go to Step~1.
\end{description}
}

The evaluation complexity analysis for the \skoffar\al{2B} algorithm is very 
closely related to that of \skoffarp, and we now discuss how the results 
of Section~\ref{complexity-s} can be adapted to the new context.

\begin{enumerate}
\item Restricting our use of the Lipschitz condition to the gradient ($p=1$), 
    Lemma \ref{lipschitz} now states that
    \beqn{Lip-f-B}
  f(x_{k+1})- \widehat \TB(x_k,\widehat s_k)
  = f(x_{k+1})- T_{f,p}(x_k,s_k) \leq \frac{\kappa_{LB}}{2} \|s_k\|^2,
  \eeqn
  and
  \[
  \|g_{k+1}-\nabla_s^1 \TB(x_k,s_k) \| \leq \kappa_{LB} \|s_k\|,
  \]
where $\kappa_{LB} \eqdef L_1+\kappa_B$.
\item Using now the decrease \req{descentB} of the model with quadratic 
regularisation, the decrease condition of Lemma~\ref{model-decrease} becomes
  \beqn{TBdecr}
  \TB(x_k,0)-\TB(x_k,s_k)
  > \frac{\sigma_k}{2} \|s_k\|^2.
  \eeqn

\item As in Lemma \ref{Ifsigmabig}, we now exploit \req{Lip-f-B} to obtain that,
if $\sigma_k \geq 2 \kappa_{LB}$, then
\beqn{sigma-upperB}
f(x_{k}) - f(x_{k+1}) > \frac{\sigma_k}{4} \|s_k\|^2.
\eeqn

\item Lemma \ref{stepgkbound} is no longer valid because its assumes 
that the regularisation order is one above that of the highest derivative 
used, while both these orders are now equal to two.  But a simple 
bound on the steplength can still be derived easily.

\llem{stepgkboundB}{
Suppose that AS.7 and AS.9 hold. At each iteration $k$, we have that
\beqn{skboundB}
\| s_k\| \leq
\frac{2\|g_k\|}{\vartheta \nu_0}.
\eeqn
}
\proof{
Using \eqref{descentB} and 
$ \widehat \mkB(\widehat s_k)= \mkB(s_k)$ it follows that
$$
\frac{1}{2}\sigma_k \|s_k\|^2\le -g_k^Ts_k-\frac{1}{2} s_k^T B_k s_k\le 
\|g_k\|\|s_k\|
$$
and the thesis follows from the fact that $\sigma_k \geq \vartheta \nu_0$.}

\item 
The proof of Lemma \ref{muk-bound} is easily adapted to the case 
where $p=1$, yielding that, for all $k \geq 0$,
\[
\mu_k \leq \max[\mu_{-1},\kappa_{LB} ].
\]
 \item The bounds \req{vkbig} and \req{k1def} may now be re-writtten as
$\nu_k \geq 2 \kappa_{LB}/\vartheta$ and
\[
k_1 \eqdef \min\left\{ k\geq 1 \mid \nu_k \geq \frac{2\kappa_{LB}}{\vartheta}\right\},
\]
respectively.

\item
Using \req{skboundB}, Lemma \ref {vk1bound} then becomes
\[
 \nu_{k_1} \leq \nu_{\max}
\eqdef 
 \frac{2\kappa_{LB}}{\vartheta}\left[1+\left(\frac{\kap{g}}{\vartheta\nu_0}\right)^2\right].
\]

\item 
The revised version of inequality \req{fmaxdef} in Lemma~\ref{fk1boundsigmak1} is now given by
\beqn{fmaxdefB}
f(x_{k_1}) \leq f_{\max} \eqdef f(x_0) +
          \frac{1}{2}\left(\frac{\kappa_{LB}}{\sigma_0} \nu_{\max}
          +\vartheta \sigma_0\right),
\eeqn
and the bound \req{sigmamaxdef} in Lemma~\ref{uppersigmak} is now
valid with
\beqn{sigmamaxdefB}
\sigma_{\max} 
\eqdef \max\left[
\bigfrac{4}{\vartheta}\left[f(x_0)-f_{\rm
  low}+\frac{1}{2}\left(\frac{\kappa_{LB}}{\sigma_0}
\nu_{\max}+\vartheta\sigma_0 \right)\right]+ \nu_{\max},
\mu_{-1}, L_1+\kappa_B, \frac{2\kappa_{LB}}{\vartheta},\nu_0 \right].
\eeqn
\item
It is of course necessary to revise our definition of a true iteration.
\begin{definition}\label{deftrueB}
Iteration $k\in \iiz{N_1(\epsilon)-1}$ is $\omega$-true whenever, 
\beqn{crucialB}
\|s_k\| \geq \omega \epsilon.
\eeqn
\end{definition}
We say
that, for some given ``preservation parameter'' $\alpha_S \in (0,1)$ and
a constant $S_{\max}>0$,
iteration $k$ is $(\alpha_S,S_{\max})$-embedded whenever, 
\beqn{true-1-B}
\|S_k g_k\|
\geq \alpha_S \|g_k\|
\tim{ and }
\|S_k\|\le S_{\max}.
\eeqn

\item
Lemma \ref{sigmanotbig} remains valid  with 
\beqn{kstardefB}
k_*
\eqdef
\left\lceil\frac{2\kappa_{LB}\epsilon^{-2}}{\vartheta \nu_0\, \omega^{2}}
\right\rceil
\tim{ and }
\sigma_k \geq 2\kappa_{LB}, \tim{ for all } k\geq k_1 
\eeqn
while Lemmas~\ref{LemmaShao} and \ref{Nits} are unchanged.
\item 
Since, for algorithm \skoffar\al{2B}, $\|g_k\| > \epsilon$ for all 
$k \leq N_1(\epsilon)-1$ (instead of $\|g_{k+1}\| > \epsilon$ for 
$k \leq N_1(\epsilon)-2$ for \skoffarp), we may continue 
to use the proof of Lemma~\ref{complexity} and obtain the following 
evaluation complexity result for the \skoffar\al{2B} algorithm.

\lthm{complexityB}{
Suppose that AS.2, AS.4, AS.6, AS.7, AS.8 and AS.9 hold, 
that $\delta_1\in (0,1)$ is given and that the \skoffar\al{2B} 
algorithm is applied to problem \req{problem}.
Define
\beqn{kapoff1B}
\kap{SKOFFAR2B}
\eqdef \frac{4\left[L_1 +\kappa_B +2(f_{\max}-f_{\rm low}) \right]}
            {\vartheta\nu_0\omega^{2}(1-\delta_1)\pione}
\eeqn
where $f_{\max}$ is defined in \req{fmaxdefB}.
Then
\[
\prob\Big[ N_1(\epsilon) \leq \kap{SKOFFAR2B} \, \epsilon^{-2}+4 \Big]
\geq \left(1 - e^{-\frac{\delta_1^2}{2}\pione k_\diamond}\right)^2
\]
where $k_\diamond=\left\lceil\frac{k_*}{(1-\delta_1)\pione}\right\rceil$
with $k_*$ given by \req{kstardefB}. 
}

Of course, using quite loose Hessian approximations in \req{modelB} has the consequence that
the complexity order is now ${\cal{O}}(\epsilon^{-2})$, which is identical to
that of other methods (such as deterministic and stochastic trust-region or 
regularisation) using the same type of approximations and objective function values.

\item 
We finally consider how Lemma \ref{useful1} can be adapted for the use of 
Gaussian scaled matrices within the \skoffar\al{2B} algorithm.

\llem{useful1B}{
Suppose  that AS.4, AS.7, AS.8 and AS.9 hold and that iteration $k\geq 0$ 
of the \skoffar\al{2B} algorithm is $(\alpha_S,S_{\max})$-embedded (in the sense of \req{true-1-B}). Then
  \[
  \|s_k\|\geq \frac{\alpha_S}{S_{\max}(\kappa_B+\theta\sigma_k)} \| g_{k}\|.
  \]
  Thus iteration $k\in\iiz{N_1(\epsilon)-1}$ is $\omega$-true (in the sense of \req{crucialB}) 
  with $\omega = \frac{\alpha_S}{S_{\max}(\kappa_B + \theta\sigma_{\max}) } $.
}

\proof{
Since
$$
S_k g_k=\nabla_s^1 \widehat \TB(x_k,\widehat s_k)-S_k B_k S_k^T \widehat s_k,
$$
we obtain from \req{gradstepB} and the definition of $(\alpha_S,S_{\max})$-embedded iteration that
\begin{eqnarray*}
    \alpha_S\|g_k\|\le \|S_k g_k\|\le S_{\max}(\kappa_B+\theta\sigma_k)\|s_k\|
\end{eqnarray*}
Using the bound \req{sigmamaxdefB} yields the desired result.
}

\noindent
We see that the constant \req{kapoff1B} now involves $S_{\max}^2$. 
In the case of scaled Gaussian matrices, \req{Gaussian-Smax} then gives a dependence 
of the constant \req{kapoff1B} in $n/\ell$, as is the case for the (non-OFFO) trust-region 
method of \cite{CartFowkShao22}.

We conclude this discussion by noting that, should the Gauss-Newton
method for nonlinear least-squares be considered, AS.3 (for $p=1$) and
AS.4 can be replaced by assuming the Lipschitz continuity of the
problem's Jacobian and the boundedness of the Jacobian and residual
(see \cite[page~295]{NoceWrig99} for a proof that this is sufficient to
ensure Lipschitz continuity and boundedness of the objective function's gradient).
\end{enumerate}

\end{document}